\documentclass[12pt,a4paper, twoside,reqno]{amsart}
\usepackage{tikz}
\usepackage{amsmath}
\usepackage{amssymb,amsfonts,mathrsfs}
\usepackage{amscd,amssymb,verbatim}
\usepackage[colorlinks=true,linkcolor=blue,citecolor=blue]{hyperref}
\usepackage{mllocal}

\theoremstyle{plain}

\numberwithin{equation}{section}

\begin{document}

\title[New proof of a  Bismut and Zhang formula]
{A new proof of a Bismut-Zhang formula for some class of representations}

\author{Maxim Braverman${}^\dag$}
\address{Department of Mathematics,
Northeastern University,
Boston, MA 02115,
USA}

\email{maximbraverman@neu.edu}
\urladdr{www.math.neu.edu/~braverman/}

\author{Boris Vertman}
\address{Mathematisches Institut,
Universit\"at Bonn,
53115 Bonn,
Germany}
\email{vertman@math.uni-bonn.de}
\urladdr{www.math.uni-bonn.de/people/vertman}

\thanks{${}^\dag$Supported in part by the NSF grant DMS-1005888.}


\begin{abstract}
Bismut and Zhang computed the ratio of the Ray-Singer and the combinatorial torsions corresponding to non-unitary representations of the fundamental group. In this note we show that for representations which belong to a connected component  containing a unitary representation the Bismut-Zhang formula follows rather easily from the Cheeger-M\"uller theorem, i.e. from the equality of the two torsions on the set of unitary representations. The proof uses the fact that the refined analytic torsion is a holomorphic function on the space of representations. 
\end{abstract}

\maketitle


\section{Introduction}\label{intro}

Let $M$ be a closed oriented odd-dimensional manifold and let $\Rep$ denote the space of representations of the fundamental group $\pi_1(M)$ of $M$. For each $\alp\in \Rep$, let $(\Ea,\na)$ be a flat vector bundle over $M$, whose monodromy representation is equal to $\alp$. We denote by $H^\b(M,\Ea)$ the cohomology of $M$ with coefficients in $\Ea$. Let $\Det(H^\b(M,\Ea))$ denote the determinant line of $H^\b(M,\Ea)$. 

Reidemeister \cite{Reidemeister35} and Franz \cite{Franz35} used a cell decomposition of $M$ to construct a combinatorial invariant of the representation $\alp\in\Rep$, called the {\em Reidemeister torsion}. In modern language it is a metric on the determinant line  $\Det(H^\b(M,\Ea))$, cf. \cite{Quillen85,BisZh92}. If $\alp$ is unitary, then this metric is independent of the cell decomposition and other choices. In general to define the Reidemeister metric one needs to make some choices. One of such choices is a Morse function $F:M\to \RR$. Bismut and Zhang \cite{BisZh92} call the metric obtain using the Morse function $F$ the {\em Milnor metric} and denote it by  $\|\cdot\|^{\M}_F$. 

Ray and Singer \cite{RaySinger71} used the de Rham complex to give a different construction of a metric on  $\Det(H^\b(M,\Ea))$. This metric is called the {\em Ray-Singer metric} and is denoted by $\|\cdot\|^{\RS}$.  Ray and Singer conjectured that the Ray-Singer and the Milnor metrics coincide  for unitary representation of the fundamental group. This conjecture was proven by Cheeger \cite{Cheeger79} and M\"uller \cite{Muller78} and extended by M\"uller \cite{Muller93} to unimodular representations. For non-unitary representations the two metrics are not equal in general. In the seminal paper \cite{BisZh92} Bismut and Zhang computed the ratio of the two metrics using very non-trivial analytic arguments.

In this note  we show that for a large class of representations  the Bismut-Zhang formula follows quite easily from the original Ray-Singer conjecture. More precisely, let $\alp_0\in\Rep$ be a unitary representation which is a regular point of the complex analytic set $\Rep$ and let $\calC\subset \Rep$ denote the connected component of $\Rep$ which contains $\alp_0$. We  derive the Bismut-Zhang formula  for all representations in $\calC$   from the Cheeger-M\"uller theorem. In other words,  we show that knowing that the Milnor and the Ray-Singer metrics coincide on unitary representations one can derive the formula for the ratio of those metrics for all representations in the  connected component  $\calC$.

The proof uses the properties of the refined analytic torsion $\rat(\alp)$ introduced in \cite{BrKappelerRATshort,BrKappelerRAT,BrKappelerRATdetline} and of the refined combinatorial torsion $\rtt(\alp)$ introduced in \cite{Turaev01,FarberTuraev00}. Both refined torsions are non-vanishing elements of the determinant line  $\Det(H^\b(M,\Ea))$ which depend holomorphically on $\alp\in \Rep$. The ratio of these sections is a holomorphic function 
\[
	\alp \ \mapsto \ \frac{\rat(\alp)}{\rtt(\alp)}
\]
on $\Rep$. We first use the Cheeger-M\"uller theorem to compute this function for unitary $\alp$.  Let now $\calC\subset \Rep$ be a connected component and suppose that a unitary representation $\alp_0$ is a regular point of $\calC$. The set of unitary representations can be viewed as the real locus of the connected complex analytic set $\calC$. As we know $\frac{\rat(\alp)}{\rtt(\alp)}$ for all points of the real locus, we can compute it for all $\alp\in \calC$ by analytic continuation. Since the Ray-Singer norm of $\rtt$ and the Milnor norm of $\rat$ are easy to compute, we obtain the Bismut-Zhang formula for all $\alp\in \calC$. 

The paper is organized as follows. In \refs{ideaBZtheorem}, we briefly outline the main steps of the proof. In \refss{Milnormetric}  and \refs{FarberTuraev torsion} we recall the construction and some properties of the Milnor metric and of the Farber-Turaev torsion. In \refs{RaySinger(rat)} we recall some properties of the refined analytic torsion. In \refs{detlinebundle} we recall the construction of the holomorphic structure on the determinant line bundle and show that the ratio of the refined analytic and the Farber-Turaev torsions is a holomorphic function on $\Rep$. Finally, in \refs{BZtheorem} we present our new proof of the Bismut-Zhang theorem for representations in the connected component $\calC$.

\section{The idea of the proof}\label{S:ideaBZtheorem}
Our proof of  the Bismut-Zhang theorem for representations in the connected component $\calC$ consists of several steps.  In this section we briefly outline these steps.

\subsection*{Step 1}
In \cite{Turaev86,Turaev90}, Turaev constructed a refined version of the combinatorial torsion associated to an acyclic representation $\alp$.  Turaev's construction depends on additional combinatorial data, denoted by $\eps$ and called the {\em Euler structure}, as well as on the {\em cohomological orientation} of $M$, i.e., on the orientation $\gro$ of the determinant line of the cohomology $H^\b(M,\RR)$ of $M$. In \cite{FarberTuraev00}, Farber and Turaev extended the definition of the Turaev torsion to non-acyclic representations. The Farber-Turaev torsion associated to a representation $\alp$, an Euler structure $\eps$, and a cohomological orientation $\gro$ is a non-zero element $\rtt(\alp)$ of the determinant line $\Det\big(H^\b(M,\Ea)\big)$. 

Let us fix a Hermitian metric $h^{\Ea}$ on $\Ea$. This scalar product induces a norm \mbox{$\|\cdot\|^\RS$} on $ \Det\big(H^\b(M,E_\alp)\big)$, called the  {\em  Ray-Singer metric}. In \refss{Milnormetric}   we use the Cheeger-M\"uller theorem to show that for unitary $\alp$
\eq{norm(rtt)}
	\|\rtt(\alp)\|^{\RS}\ = \ 1.
\end{equation}

\rem{RSFarbe-Turaev}
Theorem~10.2 of \cite{FarberTuraev00} computes the Ray-Singer norm of $\rtt$ for arbitrary representation $\alp\in \Rep$, however the proof uses the result of Bismut and Zhang,
which we want to prove here for $\alp \in \calC$.
\erem

Theorem~1.9 of \cite{BrKappelerRATdetline} computes the Ray-Singer metric of $\rat(\alp)$. Combining this result with \refe{norm(rtt)} we conclude,  cf. \refss{|rat/rtt|}, that if  $\alp$ is a unitary representation, then
\eq{RSnorm of ratio}
	\left|\frac{\rat(\alp)}{\rtt(\alp)}\right| \ = \ 
	\frac{\|\rat(\alp)\|^\RS}{\|\rtt(\alp)\|^\RS} \ =\ 1.  
\end{equation}

\subsection*{Step 2}

The Farber-Turaev torsion $\rtt(\alp)$ is a holomorphic section of the determinant line bundle 
\[
	\Dbundle\ := \  \ \bigsqcup_{\alp\in\Rep}\, \Det\big(H^\b(M,E_\alp)\big)
\]
over $\Rep$.  We denote by $\rat(\alp)/\rtt(\alp)$ the unique complex number such that
\[
	\rat(\alp)\ = \ \frac{\rat(\alp)}{\rtt(\alp)}\cdot \rtt(\alp)\ \in\ \Det\big(H^\b(M,\Ea)\big).
\]
Since both $\rtt$ and $\rat$ are holomorphic sections of $\Dbundle$, 
\[
	\alp \ \mapsto \  \frac{\rat(\alp)}{\rtt(\alp)}
\] 
is a holomorphic function on $\Rep$.

\subsection*{Step 3}
Let $\alp'$ denote the representation dual to $\alp$ with respect to a Hermitian scalar product on $\CC^n$. Then the Poincar\'e duality induces, cf. \cite[\S2.5]{FarberTuraev00} and \cite[\S10.1]{BrKappelerRATdetline}, an anti-linear isomorphism%
\footnote{There is a sign difference in the definition of the duality operator  in \cite{FarberTuraev00}  and \cite{BrKappelerRATdetline}, which is not essential for the discussion in this paper.}
\[
	D:\,\Det\big(H^\b(M,\Ea)\big) \ \longrightarrow \ \Det\big(H^\b(M,E_{\alp'})\big).
\]
In particular, when $\alp$ is a unitary representation, $D$ is an anti-linear automorphism of $\Det\big(H^\b(M,\Ea)\big)$. Hence,
\eq{Drat/Drtt0}
	\frac{D(\rat(\alp))}{D(\rtt(\alp))} \ = \ \overline{\frac{\rat(\alp)}{\rtt(\alp)}}.
\end{equation}

Using Theorem~7.2 and formula (9.4) of \cite{FarberTuraev00} we compute the ratio $D(\rtt(\alp))/\rtt(\alp)$, cf. \refe{Drtt=ceps unitary} (here $\alp$ is a unitary representation). On the analytic side Theorem~10.3 of \cite{BrKappelerRATdetline} computes the ratio $D(\rat(\alp))/\rat(\alp)$. Combining these two results we get
\eq{Drat/Drtt=rat/rtt}
	\overline{\frac{\rat(\alp)}{\rtt(\alp)}}\ = \ f_2(\alp)\cdot \frac{\rat(\alp)}{\rtt(\alp)},
\end{equation}
where $f_2$ is a function on $\Rep$ computed explicitly in \refe{Drat/Drtt}. 

From \refe{Drat/Drtt0} and \refe{Drat/Drtt=rat/rtt} we conclude that 
\eq{rat/rtt=}
	\left(\frac{\rat(\alp)}{\rtt(\alp)}\right)^2 \ = \  f_1(\alp)^2\cdot f_2(\alp)
\end{equation}
for any unitary representation $\alp$, cf. \refe{rat/rtt unitary}, where $f_1(\alp)= \rat(\alp) / \rtt(\alp)$.

\subsection*{Step 4}
The right hand side of \refe{rat/rtt=} is an explicit function of a unitary representation $\alp$. It turns out that it is a restriction of a holomorphic function $f(\alp)$ on $\Rep$ to the set of unitary representations.  Recall that the connected component  $\calC$ contains a regular point which is a unitary representation.  The set of unitary representations can be viewed as the {\em real locus} of the complex analytic set $\calC$. Hence any two holomorphic functions which coincide on the set of unitary representations, coincide on $\calC$. We conclude now from \refe{rat/rtt=} that 
\eq{rat/rtt=f}
		\left(\frac{\rat(\alp)}{\rtt(\alp)}\right)^2 \ = \ f(\alp), \qquad\text{for all}\ \ \alp\in \calC. 
\end{equation}

\subsection*{Step 5}
Recall that we denote by $\|\cdot\|^M_F$ the Milnor metric associated to the Morse function $F$. In \refs{FarberTuraev torsion} we compute  the Milnor metric 
\eq{milnor=}
	\|\rtt(\alp)\|^M_F\ = \ h_1(\alp),
\end{equation}
where  $h(\alp)$ is a real valued function on $\Rep$ given explicitly by the right hand side of \refe{MilnorofTuraev}

Theorem~1.9 of \cite{BrKappelerRATdetline} computes the Ray-Singer norm of the refined analytic torsion:
\begin{equation}\label{E:RSrat}
	\|\rat(\alp)\|^\RS \ = \ h_2(\alp),
\end{equation}
where  $h_2(\alp)$ is a real valued function on $\Rep$ given explicitly by the right hand side of \refe{|rat|}. 
Combining  \refe{rat/rtt=f}  with \refe{RSrat}, we get
\eq{BZformula1}
	\frac{\|\cdot\|^\RS}{\|\cdot\|^{\M}_F}\ = \ \frac{\|\rat(\alp)\|^\RS}{\|\rtt(\alp)\|^{\M}_F}\cdot \left|\,\frac{\rtt(\alp)}{\rat(\alp)} \,\right|
	\ = \ \frac{h_2(\alp)}{h_1(\alp)\cdot|f(\alp)|}.
\end{equation}
This is exactly the Bizmut-Zhang formula \cite[Theorem~0.2]{BisZh92}.

The rest of the paper is occupied with the details of the proof outlined above.

\section{The  Milnor metric and the Farber-Turaev torsion}\label{S:FarberTuraev torsion}

In this section we briefly recall the definitions and the main properties of the Milnor metric and the Farber-Turaev refined combinatorial torsion. We also compute the Milnor norm of the Farber-Turaev torsion. 
\subsection{The Thom-Smale complex}\label{SS:Thom-Smale}
Set
\[
	C^k(K,\Ea) \ = \ \bigoplus_{\substack{x\in Cr(F)\\ \ind_F(x)=k}} E_{\alp,x}, \qquad k=1\nek n,
\]
where $E_{\alp,x}$ denotes the fiber of $E_\alp$ over $x$ and the direct sum is over the critical points $x\in Cr(F)$ of the Morse function $F$ with Morse-index $\ind_F(x)=k$. If the Morse function is $F$ generic, then using the gradient flow of $F$ one can define the {\em Thom-Smale complex}  $(C^\b(K,\Ea),\d)$ whose cohomology is canonically isomorphic to $H^\b(M,\Ea)$, cf. for example \cite[\S{}I c]{BisZh92}.

\subsection{The Euler structure}\label{SS:Euler}
The {\em Euler structure} $\eps$ on $M$ can be described as (an equivalence class of) a pair $(F,c)$ where $c$ is a 1-chain in $M$ such that 
\eq{dc=cr}
	\d c \ =\ \sum_{x\in Cr(F)} (-1)^{\ind_F(x)}\cdot x,
\end{equation}
cf. \cite[\S3.1]{BurgheleaHaller_Euler}. We denote the set of Euler structures on $M$ by $\Eul$.

\rem{Euler structue}
The Euler structure was introduced by Turaev \cite{Turaev90}. Turaev presented several equivalent definitions and the equivalence of these definitions is a nontrivial result. Burghelea and Haller \cite{BurgheleaHaller_Euler} found a very nice way to unify these definitions. They suggested a new definition which is obviously equivalent to the two definitions of Turaev. In this paper we use the definition introduced by Burghelea and Haller. 
\erem

\subsection{The Kamber-Tondeur form}\label{SS:KamberTondeur}
To define the Milnor and the Ray-Singer metrics on $\Det(H^\b(M,\Ea))$ we fix a Hermitian metric $h^{\Ea}$ on $\Ea$. This metric is not necessary flat and the measure of non-flatness is given by taking the trace of $(h^{\Ea})^{-1}\na h^{\Ea} \in \Ome^1(M, \textup{End} E_\alp)$ which defines the {\em Kamber-Tondeur form}
\begin{equation}\label{E:KamberTondeur}
	\tet(h^{\Ea})\ := \ \Tr\,\big[\, (h^{\Ea})^{-1}\na h^{\Ea}\,\big] \ \in \ \Ome^1(M),
\end{equation}
cf. \cite{KamberTondeur67} (see also \cite[Ch.~IV]{BisZh92}). 

Let $\Det(\Ea)\to M$ denote the determinant line bundle of $\Ea$, i.e. the line bundle whose fiber over $x\in M$ is equal to the determinant line $\Det(E_{\alp,x})$ of the fiber $E_{\alp,x}$ of $\Ea$. The connection $\na$ and the metric $h^{\Ea}$ induce a flat connection $\na^{\Det}$ and a metric $h^{\Det(\Ea)}$ on $\Det(\Ea)$. Then 
\eq{tetE=tetDetE}
	\tet(h^{\Det(\Ea)}) \ = \ \tet(h^{\Ea}).
\end{equation}
For a curve $\gam:[a,b]\to M$ let 
\eq{monodromy det}
	\alp(\gam):\,E_{\alp,\gam(a)}\ \to \ E_{\alp,\gam(b)}, \quad 
	\alp^{\Det}(\gam):\,\Det(E_{\alp,\gam(a)})\ \to \ \Det(E_{\alp,\gam(b)})
\end{equation}
denote the parallel transports along $\gam$. Then 
\eq{Detalp=alpDet}
	\Det\big(\alp(\gam)\big)\ = \ \alp^{\Det}(\gam).
\end{equation}

Let $\tilgam(t)\in \Det(E_{\alp,\gam(t)})$ denote the horizontal lift of the curve $\gam$. By the definition of the Kamber-Tondeur form we have
\eq{logh/h}
	\log\, \frac{h^{\Det(\Ea)}(\tilgam(b),\tilgam(b))}{h^{\Det(\Ea)}(\tilgam(a),\tilgam(a))}
	\ = \ \int_\gam\, \tet(h^{\Det(\Ea)})\ = \ \int_\gam\,\tet(h^{\Ea}),
\end{equation}
where in the last equality we used \refe{tetE=tetDetE}. 

If $\gam$ is a closed curve, $\gam(a)=\gam(b)$, we obtain 
\[
	 \frac{h^{\Det(\Ea)}(\tilgam(b),\tilgam(b))}{h^{\Det(\Ea)}(\tilgam(a),\tilgam(a))} \ = \ 
	 \big|\alp^{\Det}(\gam)\big|^2\ = \ \big|\Det\big(\alp(\gam)\big)\big|^2.
\]
Hence from \refe{logh/h} we obtain
\eq{|DetPsi|}
	\big|\Det\big(\alp(\gam)\big)\big|\ = \ e^{\frac12\int_\gam\tet(h^{\Ea})}.
\end{equation}
\subsection{The  Milnor metric}\label{SS:Milnormetric} 
The Hermitian metric $h^{\Ea}$ on $\Ea$ defines a scalar product on the spaces $C^\b(K,\Ea)$ and, hence, a metric  $\|\cdot\|_{\Det(C^\b(K,\Ea))}$ on the determinant line of $C^\b(K,\Ea)$. Using the isomorphism 
\begin{equation}\label{E:phialp2}
    \phi:\, \Det\big(\,C^\b(K,\Ea)\,\big) \ {\longrightarrow} \ \Det\big(\,H^\b(M,\Ea)\,\big),
\end{equation}
cf. formula (2.13) of \cite{BrKappelerRATdetline}, we thus obtain a metric on $\Det\big(\,H^\b(M,\Ea)\,\big)$, called the {\em Milnor metric} associated with the Morse function $F$ and denoted by $\|\cdot\|^{\M}_F$.

\subsection{The Farber-Turaev torsion}\label{SS:FarberTuraev torsion}
Turaev \cite{Turaev90} showed  that if an Euler structure is fixed, then  the scalar product on the spaces $C^k(K,\Ea)$  allows one to construct  not only a metric on the determinant line $\Det\big(C^\b(K,\Ea)\big)$ but also an element of this line, defined modulo sign. 

We recall briefly Turaev's construction.  Fix a base point $x_*\in M$. Then every Euler structure $\eps$ can be represented by a pair $(F,c)$ such that 
\[
	c\ = \ \sum_{x\in Cr(F)}\, (-1)^{\ind_F(x)}\gam_x,
\]
with $\gam_x:[0,1]\to M$ being a smooth curve such that $\gam_x(0)= x_*$ and $\gam_x(1)= x$. The chain $c$ is often referred to as a {\em Turaev spider}.

We need to construct an element of the the determinant line  $\Det\big(C^\b(K,\Ea)\big)$  of the cochain complex $C^\b(K,\Ea)$. It is easier to start with constructing an element in the determinant line of the {\em chain} complex. Since the cochain complex is dual to the chain complex of the bundle $E_{\alp'}$, where  $\alp'$ denote the representation dual to $\alp$, we construct an element in the  determinant line $\Det\big(C_\b(K,E_{\alp'})\big)$. This is done as follows: 

Fix an element $v_*\in \Det(E_{\alp',x_*})$ whose norm with respect to the Hermitian metric $h^{\Det(E_{\alp'})}$ is equal to 1 and set 
\[
	v_x\ := \ {\alp'}^{\Det}(\gam_x)(v_*)\ \in \ \Det(E_{\alp',x}),
\]
where ${\alp'}^{\Det}$ is the monodromy of the induced connection on the determinant line bundle $\Det(E_{\alp'})$, cf. \refe{monodromy det}. Let 
\[
	|v|^{\Det(E_{\alp'})} \ := \ \sqrt{h^{\Det(E_{\alp'})}(v,v)}
\]
denote the norm induced on $\Det(E_{\alp'})$ by the Hermitian metric $h^{\Det(E_{\alp'})}$. Then from  \refe{logh/h} we obtain
\eq{h(v_x)}
	|v_x|^{\Det(E_{\alp'})}  \ = \ \frac{|v_x|^{\Det(E_{\alp'})}}{|v_*|^{\Det(E_{\alp'})}}
	\ = \  e^{\frac12\int_{\gam_x}\tet( h^{\Det(E_{\alp'})})}
	\ =\ e^{-\frac12\int_{\gam_x}\tet( h^{\Det(E_{\alp})} )}.
\end{equation}

Let
\[
	v\ = \ \prod_{x\in Cr(F)} v_x^{(-1)^{\ind_F(x)}} \in\  \Det\big(C_\b(K,E_{\alp'})\big)/\pm.
\] 
(The sign indeterminacy comes form the choice of the order of the critical points of $F$.)
From \refe{h(v_x)} we conclude that 
\eq{norm of v}
	\|v\|_{\Det(C_\b(K,E_{\alp'}))} \ = \ e^{-\frac12\int_c\tet( h^{\Det(E_{\alp})} )}.
\end{equation}

Let $\<\cdot,\cdot\>$ denote the natural pairing
\[
	\Det\big(C^\b(K,E_{\alp})\big)\times \Det\big(C_\b(K,E_{\alp'})\big)\ \to \ \CC
\]
and let $\nu\in \Det\big(C^\b(K,E_{\alp})\big)/\pm$  be the unique element such that $\<\nu,v\>=1$. From \refe{norm of v} we now obtain 
\eq{norm of nu}
	\|\nu\|_{\Det(C^\b(K,E_{\alp}))} \ = \ e^{\frac12\int_c\tet( h^{\Det(E_{\alp})} )}.
\end{equation}

Using the isomorphism \refe{phialp2} we obtain an element 
\eq{phi(v)}
	\phi(\nu)\ \in \ \Det\big(H^\b(M,\Ea)\big)/\pm.
\end{equation}
To fix the sign one can choose a {\em cohomological orientation} $\gro$, i.e. an orientation of the determinant line $\Det(H^\b(M,\RR)$. Thus, given the Euler structure $\eps$ and the cohomological orientation $\gro$ we obtain a sign refined version of $\phi(\nu)$  which we call the {\em Farber-Turaev torsion} and denote by 
\eq{rhoFT}
	\rtt(\alp) \in \  \Det\big(H^\b(M,\Ea)\big).
\end{equation}

\subsection{The Milnor norm of the Farber-Turaev torsion}\label{SS:Milnor Farber-Turaev}
From \refe{norm of nu} we immediately get
\eq{MilnorofTuraev}
	\|\rtt(\alp)\|_{F}^{\M} \ = \ e^{\frac12\int_c\tet(h^{\Ea})}.
\end{equation}
In particular, if $\alp$ is a unitary representation, then $h^{\Ea}$ is a flat Hermitian metric and $\tet(h^{\Ea})=0$. Hence, if $\alp$ is unitary, then
\eq{MilnorofTuraevunit}
	\|\rtt(\alp)\|_{F}^{\M} \ = \ 1.
\end{equation}
We now use the Cheeger-M\"uller theorem to conclude that 
\eq{RaySingerofTuraevunit}
	\|\rtt(\alp)\|^{\RS} \ = \ 1, \qquad \text{if $\alp$ is unitary}.
\end{equation}

\subsection{Dependence of the Farber-Turaev torsion on the Euler structure}\label{rtt depends on e}
For a homology class $h\in H_1(M,\ZZ)$ and an Euler structure $\eps=(F,c)\in \Eul$ we set 
\eq{H1 acts on Eul}
	h\eps\ := \ (F,c+h)\ \in \ \Eul.
\end{equation}
This defines a free and transitive action of $H_1(M,\ZZ)$ on $\Eul$,  cf.  \cite[\S5]{FarberTuraev00} or \cite[\S3.1]{BurgheleaHaller_Euler}.

One easily checks, cf. \cite[page~211]{FarberTuraev00},  that 
\eq{rtt he}
	\rho_{h\eps,\gro}(\alp)\ = \ \Det \big(\alp(h)\big)\cdot \rtt(\alp).
\end{equation}
From \refe{|DetPsi|} and \refe{MilnorofTuraev} we now obtain
\eq{Milnorrhohe}
	\|\rho_{h\eps,\gro}(\alp)\|^{\M}_F\ = \ e^{-\frac12\int_{c+h}\tet(h^{\Ea})}.
\end{equation}

\section{The Ray-Singer norm of the Refined Analytic Torsion}\label{S:RaySinger(rat)}

In \cite{BrKappelerRATdetline} Braverman and Kappeler defined an element of $\Det(H^\b(M,\Ea))$ called the {\em refined analytic torsion} and denoted by $\rat(\alp)$. They also computed the Ray-Singer norm $\|\rat(\alp)\|^{\RS}$ of the refined analytic torsion. In this section we recall the result of this computation.

\subsection{The odd signature operator}\label{SS:oddsign}
Fix a Riemannian metric $g^M$ on $M$ and let $*:\Ome^\b(M,\Ea)\to \Ome^{m-\b}(M,\Ea)$ denote the Hodge $*$-operator, where  $m=\dim M$. Define the {\em chirality
operator} 
\[
	\Gam\ =\  \Gam(g^M):\,\Ome^\b(M,\Ea)\ \to\ \Ome^\b(M,\Ea)
\] 
by the formula
\begin{equation}\label{E:Gam}
    \Gam\, \ome \ := \ i^r\,(-1)^{\frac{k(k+1)}2}\,*\,\ome, \qquad \ome\in \Ome^k(M,E),
\end{equation}
where\/ $r=\frac{m+1}2$. The numerical factor in \refe{Gam} has been chosen so that $\Gam^2=1$, cf. Proposition~3.58 of \cite{BeGeVe}.

The {\em odd signature operator} is the operator
\begin{equation} \label{E:oddsignGam}
    \B\ = \  \B(\na,g^M) \ := \ \Gam\,\na \ + \ \na\,\Gam:\,\Ome^\b(M,\Ea)\ \longrightarrow \  \Ome^\b(M,\Ea).
\end{equation}

\subsection{The eta invariant}\label{SS:eta}
We recall from \cite[\S3]{BrKappelerRATdetline} the definition of the sign-refined $\eta$-invariant $\eta(\na,g^M)$ of the (not necessarily unitary) connection $\na$. 

 Let $\Pi_>$  (resp. $\Pi_<$) be the projection whose image contains the span of all
generalized eigenvectors of $\B$ corresponding to eigenvalues $\lam$ with $\RE\lam>0$ (resp. with $\RE\lam<0$) and whose kernel contains the span of
all generalized eigenvectors of $\B$ corresponding to eigenvalues $\lam$ with $\RE\lam\le0$ (resp. with $\RE\lam\ge0$), cf.
\cite[Appendix~B]{Ponge-asymetry}. We define the $\eta$-function of $\B$ by the formula
\eq{eta}
    \eta_{\tet}(s,\B) \ = \ \Tr\,\big[ \Pi_>\B^s_\tet\big] \ - \ \Tr\,\big[ \Pi_<(-\B)^s_\tet\big], 
\end{equation}
where $\tet$ is an Agmon angle for both operators $\B$ and $-\B$ and $\B^s_\tet$ denotes the complex power of $\B$ defined relative to the spectral cut along the ray $\{re^{i\tet}:r>0\}$, cf. \cite{Seeley67,ShubinPDObook}.
It was shown by Gilkey, \cite{Gilkey84}, that $\eta_\tet(s,\B)$ has a meromorphic extension to the whole complex plane $\CC$ with isolated simple
poles, and that it is regular at $s=0$. Moreover, the number $\eta_\tet(0,\B)$ is independent of the Agmon angle $\tet$.

Let $m_+(\B)$ (resp., $m_-(\B)$) denote the number of eigenvalues of $\B$, counted with their algebraic multiplicities, on the positive (resp.,
negative) part of the imaginary axis. Let $m_0(\B)$ denote algebraic multiplicity of 0 as an eigenvalue of $\B$.

\defe{etainv}
The $\eta$-invariant $\eta(\na,g^M)$ of the pair $(\na,g^M)$ is defined by the formula
\eq{etainv}
    \eta(\na,g^M) \ = \ \frac{\eta_\tet(0,\B)+m_+(\B)-m_-(\B)+m_0(\B)}2.
\end{equation}
\edefe
If the representation $\alp$ is unitary, then the operator $\B$ is self-adjoint and $\eta(\na,g^M)$ is real. If $\alp$ is not unitary then, in general, $\eta(\na,g^M)$ is a complex number. Notice, however, that while the real part of $\eta(\na,g^M)$ is a non-local spectral invariant, the imaginary part $\IM{}\eta(\na,g^M)$ of $\eta(\na,g^M)$ is local and relatively easy to compute, cf. \cite{Gilkey84,MaZhang06eta}.

We also note that the imaginary part  of the $\eta$-invariant is independent of the Riemannian metric $g^M$.

\subsection{The Ray-Singer norm of the refine analytic torsion}\label{SS:|rat|}
Let  $\eta(\na,g^M)$ denote the $\eta$-invariant of the odd signature operator corresponding to the connection $\na$. By Theorem~1.9 of \cite{BrKappelerRATdetline} 
\eq{|rat|}
	\|\rat(\alp)\|^{\RS}\ = \ e^{\pi\IM(\eta(\na,g^M))}.
\end{equation}
In particular, when $\alp$ is a unitary representation,  $\eta(\na,g^M)$ is real and we get
\eq{|rat| unitary}
	\|\rat(\alp)\|^{\RS}\ = \ 1.
\end{equation}

\section{The Determinant Line Bundle over the Space of Representations}\label{S:detlinebundle}

The space $\Rep$ of complex $n$-dimensional representations of $\p$ has a natural structure of a complex analytic space, cf., for example,
\cite[\S13.6]{BrKappelerRAT}. The disjoint union
\begin{equation}\label{E:Detbundle}
    \Dbundle \ := \ \bigsqcup_{\alp\in\Rep}\, \Det\big(H^\b(M,E)\big)
\end{equation}
is a line bundle over $\Rep$, called the {\em determinant line bundle}. In \cite[\S3]{BrKappelerRATdetline_hol}, Braverman and Kappeler constructed a natural holomorphic structure on $\Dbundle$,  with respect to which both the refined analytic torsion $\rat(\alp)$ and the Farber-Tureav torsion $\rtt(\alp)$ are holomorphic sections. In this section we first recall this construction and then consider the ratio $\rat/\rtt$ of these two sections as a holomorphic function on $\Rep$.

\subsection{The flat vector bundle induced by a representation}\label{SS:Ealp}
Denote by $\pi:\tilM\to M$ the universal cover of $M$. For $\alp\in \Rep$, we denote by
\begin{equation}\label{E:Ealp}
  E_\alp \ := \ \tilM\times_\alp \CC^n \ \longrightarrow M
\end{equation}
the flat vector bundle induced by $\alp$. Let $\n_\alp$ be the flat connection on $E_\alp$ induced from the trivial connection on
$\tilM\times\CC^n$. 

For each connected component (in classical topology) $\calC$ of $\Rep$, all the bundles $E_\alp$, $\alp\in \calC$, are isomorphic, see e.g.
\cite{GoldmanMillson88}.

\subsection{The combinatorial cochain complex}\label{SS:CKalp}
Fix a CW-decomposition $K=\{e_1\nek e_N\}$ of $M$. For each $j=1\nek N$, fix a lift $\tile_j$, i.e., a cell of the CW-decomposition of $\tilM$,
such that $\pi(\tile_j)= e_j$. By \refe{Ealp}, the pull-back of the bundle $E_\alp$ to $\tilM$ is the trivial bundle $\tilM\times\CC^n\to
\tilM$. Hence, the choice of the cells $\tile_1\nek \tile_N$ identifies the cochain complex $C^\b(K,\alp)$ of the CW-complex $K$ with
coefficients in $\Ea$ with the complex
\begin{equation}\label{E:C(K,alp)}
    \begin{CD}
        0 \ \to \CC^{n\cdot k_0} @>{\pa_0(\alp)}>> \CC^{n\cdot k_1} @>{\pa_1(\alp)}>>\cdots @>{\pa_{m-1}(\alp)}>> \CC^{n\cdot k_m} \ \to \ 0,
    \end{CD}
\end{equation}
where $k_j\in \ZZ_{\ge0}$  ($j=0\nek m$) is equal to the number of $j$-dimensional cells of $K$ and the differentials $\pa_j(\alp)$ are $(nk_{j}\times{}nk_{j-1})$-matrices depending analytically on $\alp\in \Rep$.

The cohomology of the complex \refe{C(K,alp)} is canonically isomorphic to $H^\b(M,\Ea)$. Let
\begin{equation}\label{E:phialp}
    \phi_{C^\b(K,\alp)}:\, \Det\big(\,C^\b(K,\alp)\,\big) \ {\longrightarrow} \ \Det\big(\,H^\b(M,\Ea)\,\big)
\end{equation}
denote the natural  isomorphism between the determinant line of the complex and the determinant line of its cohomology, cf.  \cite[\S2.4]{BrKappelerRATdetline}
\subsection{The holomorphic structure on $\Dbundle$}\label{SS:holonDet}
The standard bases of $\CC^{n\cdot k_j}$ ($j=0\nek m$) define an element $c\in \Det\big(C^\b(K,\alp)\big)$, and, hence, an isomorphism
\[
    \psi_\alp:\,\CC \ {\longrightarrow} \ \Det\big(C^\b(K,\alp)\big), \qquad z\ \mapsto z\cdot c.
\]
Then the map
\begin{equation}\label{E:combsection}
    \sig:\, \alp \ \mapsto \ \phi_{C^\b(K,\alp)}\big(\,\psi_\alp(1)\,\big) \ \in \ \Det\big(\,H^\b(M,\Ea)\,\big),
\end{equation}
where $\alp\in \Rep$ is a nowhere vanishing section of the determinant line bundle $\Dbundle$ over $\Rep$.
\begin{Def}\label{D:holsection}
We say that a section $s(\alp)$ of\/ $\Dbundle$ is {\em holomorphic} if there exists a holomorphic function $f(\alp)$ on $\Rep$, such that
$s(\alp)= f(\alp)\cdot{}\sig(\alp)$.
\end{Def}
This defines a holomorphic structure on $\Dbundle$, which is independent of the choice of the lifts $\tile_1\nek \tile_N$ of $e_1\nek e_N$,
since for a different choice of lifts the section $\sig(\alp)$ will be multiplied by a constant. It is shown in \cite[\S3.5]{BrKappelerRATdetline_hol} that this holomorphic structure is also independent of the CW-decomposition $K$ of $M$.

\th{torsion is holomorphic}
Both the refined analytic torsion $\rat(\alp)$ and the Farber-Turaev torsion $\rtt(\alp)$ are holomorphic sections of\, $\Dbundle$ with respect to the holomorphic structure described above. 
\eth
\prf
The fact that the Farber-Turaev torsion is holomorphic is established in  Proposition~3.7 of \cite{BrKappelerRATdetline_hol}. The fact that the refined analytic torsion is holomorphic is proven in Theorem~4.1 of \cite{BrKappelerRATdetline_hol}.
\eprf

\subsection{The ratio of the torsions as a holomorphic function}\label{SS:rat/rtt} 
Since both $\rtt$ and  $\rat$ are holomorphic nowhere vanishing section of  the same line bundle there exists a holomorphic function 
\[
	\kap:\,\Rep\ \to\ \CC\backslash\{0\}
\]
such that 
\[
	\rat(\alp)\ = \ \kap(\alp)\cdot \rtt(\alp).
\]
We shall denote this function by
\eq{rat/rtt}
	\kap(\alp)=\frac{\rat(\alp)}{\rtt(\alp)}.
\end{equation}

\subsection{The absolute value of $\frac{\rat(\alp)}{\rtt(\alp)}$ for unitary representations}\label{SS:|rat/rtt|}
Combining \refe{|rat| unitary} with \refe{RaySingerofTuraevunit} we obtain 
\eq{RSnorm of ratio2}
	\left|\frac{\rat(\alp)}{\rtt(\alp)}\right| \ = \ 
	\frac{\|\rat(\alp)\|^\RS}{\|\rtt(\alp)\|^\RS} \ =\ 1, \qquad \text{if $\alp$ is unitary}.
\end{equation}

\section{The Bismut-Zhang theorem for some non-unitary representations}\label{S:BZtheorem}

We now present our proof of the Bismut-Zhang theorem \cite[Theorem~0.2]{BisZh92} for representations in the connected component $\calC$.

\subsection{The duality operator}\label{SS:duality}
Let $\alp'$ denotes the representation dual to $\alp$. The Poincar\'e duality defines a non-degenerate pairing 
\[
	\Det\big(H^k(M,\Ea\big)\times \Det\big(H^{m-k}(M,E_{\alp'})\big)\ \to \CC, \qquad k=0\nek m,
\] 
and, hence, an anti-linear map 
\eq{duality}
	D: \Det\big(H^\b(M,\Ea)\big) \ \to \ \Det\big(H^\b(M,E_{\alp'})\big)
\end{equation}
see \cite[\S2.5]{FarberTuraev00} and \cite[\S10.1]{BrKappelerRATdetline} for details. 

By Theorem~10.3 of \cite{BrKappelerRATdetline}  we have 
\eq{Drat=rat}
	    D\,\rat(\alp) \ = \ \rat(\alp')\cdot e^{ 2i\pi\,\big(\eta(\na,g^M) - (\rank E)\,\eta_\trivial(g^M)\big)},
\end{equation}
where $\eta(\na,g^M)$ is defined in \refd{etainv} and $\eta_\trivial$ is the $\eta$-invariant corresponding to the standard connection on the trivial line bundle $M\times\CC\to M$.

\subsection{The dual of the Farber-Turaev torsion}\label{SS:Drtt}
By Theorem~7.2  of \cite{FarberTuraev00}
\eq{Drtt=rtt}
	    D\,\rtt(\alp) \ = \ \pm\,\rho_{\eps^*,\gro}(\alp'),
\end{equation}
where $\eps^*:=(-F,-c)$ is the {\em dual Euler structure} on $M$. 

We shall use formula \refe{rtt he} in the following situation: if $\eps = (F,c)\in \Eul$ then the Euler structure $\eps^*:=(-F,-c)$ is called {\em dual} to $\eps$.  Since $H_1(M,\ZZ)$ acts freely and transitively on $\Eul$ there exists $c_\eps\in H_1(M,\ZZ)$ such that 
\eq{eps=ceps*}
	\eps\ = \ c_\eps\,\eps^*.
\end{equation}
The homology class $c_\eps$ was introduced by Turaev \cite{Turaev90} and is called the {\em characteristic class of the Euler structure}. From \refe{rtt he}  and \refe{Drtt=rtt} we now conclude that 
\eq{Drtt=ceps}
	    D\,\rtt(\alp) \ = \ \pm\rho_{\eps^*,\gro}(\alp')\ = \ \pm\Det \big(\alp'(c_\eps)\big)\cdot\rho_{\eps,\gro}(\alp').
\end{equation}
If $\alp$ is a unitary representation, then $\alp=\alp'$. Hence, it follows from \refe{Drtt=ceps} that 
\eq{Drtt=ceps unitary}
	   \rho_{\eps,\gro}(\alp') \ = \  \pm\big(\Det \big(\alp(c_\eps)\big)\big)^{-1}\cdot D\,\rtt(\alp).
\end{equation}

\subsection{The ratio of torsions for unitary representations}\label{SS:ratiounitary}
Combining \refe{Drat=rat} and \refe{Drtt=ceps unitary} we conclude that for unitary $\alp$
\eq{Drat/Drtt}
	\frac{D\rat(\alp)}{D\rtt(\alp)}\ = \  
	\pm\Det \big(\alp(c_\eps)\big)\cdot e^{ 2i\pi\,\big(\eta(\na,g^M) - (\rank E)\,\eta_\trivial(g^M)\big)}\cdot
	\frac{\rat(\alp)}{\rtt(\alp)}.
\end{equation}
Since $D$ is an anti-linear involution we have 
\[
	\frac{D\rat(\alp)}{D\rtt(\alp)}\ = \  \overline{\frac{\rat(\alp)}{\rtt(\alp)}}.
\]
Hence, it follows from \refe{Drat/Drtt} that 
\eq{rat/rtt= |rat/rtt|}
	\left(\frac{\rat(\alp)}{\rtt(\alp)}\right)^2\ = \ \pm\Det \big(\alp(c_\eps)\big)\cdot e^{ -2i\pi\,\big(\eta(\na,g^M) - (\rank E)\,\eta_\trivial(g^M)\big)}\cdot
	\left|\frac{\rat(\alp)}{\rtt(\alp)}\right|^2.
\end{equation}
Combining this equality with \refe{RSnorm of ratio2} we obtain for unitary $\alp$
\eq{rat/rtt unitary}
 \left(\frac{\rat(\alp)}{\rtt(\alp)}\right)^2\ = \ \pm\Det \big(\alp(c_\eps)\cdot e^{ -2i\pi\,\big(\eta(\na,g^M) - (\rank E)\,\eta_\trivial(g^M)\big)}.
\end{equation}

\subsection{The ratio of torsions for non-unitary representations}\label{SS:rationonunitary}
Suppose now that  $\calC\subset \Rep$ is a connected component and $\alp_0\subset \calC$ is a unitary representation which is a regular point of the complex analytic set $\calC$. The set of unitary representations is the fixed point set of the anti-holomorphic involution 
\[
	\tau:\,\Rep\ \to \Rep, \qquad \tau:\,\alp\ \mapsto \ \alp'.
\]
Hence it is a totally real submanifold of $\Rep$ whose real dimension is equal to $\dim_\CC\calC$, see for example \cite[Proposition~3]{Ho04}. In particular there is a holomorphic coordinates system $(z_1\nek z_r)$ near $\alp_0$ such that the unitary representations form a {\em real neighborhood} of $\alp_0$, i.e. the set $\IM{z_1}=\ldots =\IM{z_r}=0$. Therefore, cf. \cite[p.~21]{ShabatSeveralVariables}, if two holomorphic functions coincide on the set of unitary representations they also coincide on $\calC$. We conclude that the equation \refe{rat/rtt unitary} holds for all representations $\alp\in \calC$. Hence,  using \refe{|rat|} and \refe{MilnorofTuraev} we obtain for every $\alp\in \calC$
\eq{ration of metrics}
	\frac{\|\cdot\|^\RS}{\|\cdot\|^{\M}_F}\ = \ \frac{\|\rat(\alp)\|^\RS}{\|\rtt(\alp)\|^{\M}_F}\cdot \left|\,\frac{\rtt(\alp)}{\rat(\alp)} \,\right|
	\ = \ 
	\big|\Det \big(\alp(c_\eps)\big)\big|^{-1/2}\cdot e^{-\frac12\int_c\tet(h^{\Ea})}.
\end{equation}

\subsection{The absolute value of the determinant of $\alp(c_\eps)$}\label{SS:detseps}
Let 
\[
	\PD:\,H_1(M,\RR)\ \to\ H^{n-1}(M,\RR)
\] 
denote the Poincar\'e isomorphism. By Proposition~3.9 of \cite{BurgheleaHaller_Euler} there exists a map
\[
	P:\Eul \ \to \ \Ome^{n-1}(M,\RR)	
\]
such that 
\eq{P(heps)}
\begin{aligned}
	P(h\eps)\ &= \ P(\eps)\ +\ \PD(h),\\
	P(\eps^*)\ &= \ -P(\eps),
\end{aligned}
\end{equation}
and if $\eps=(X,c)$ then for every $\ome\in \Ome^1(M,\RR)$
\eq{intc=}
	\int_c\,\ome \ = \  \int_M\,\ome\wedge X^*\Psi(g)\ - \ \int_M\,\ome\wedge P(\eps).
\end{equation}
Here $\Psi(g)$ is the Mathai-Quillen current on $TM$, cf. \cite[\S{}III c]{BisZh92} and $X^*\Psi(g)$ denotes the pull-back of this current by $X:M\to TM$.

Combining \refe{eps=ceps*} with \refe{P(heps)} we obtain
\[
	P(\eps)\ = \ P(\eps^*)\ + \ \PD(c_\eps)\ = \ -P(\eps) \ + \ \PD(c_\eps).
\]
Thus
\eq{ceps=}
	\PD(c_\eps)\ = \ 2\,P(\eps).
\end{equation}
Combining this equality with \refe{intc=} we get
\eq{inttc=}
	\int_c\,\ome \ = \  \int_M\,\ome\wedge X^*\Psi(g)\ - \ \frac12\,\int_M\,\ome\wedge \PD(c_\eps).
\end{equation}
Notice now that 
\[
	\int_M\,\ome\wedge \PD(c_\eps) \ = \ \int_{c_\eps}\,\ome.
\]
Hence, from  \refe{inttc=} we obtain
\eq{intceps=}
	\int_{c_\eps}\,\ome\ =\ -2\,\int_c\,\ome \ + \  2\,\int_M\,\ome\wedge X^*\Psi(g).
\end{equation}
In particular, setting $\ome=\tet(h^{\Ea})$ and using \refe{|DetPsi|} we obtain
\eq{detalpceps}
	\big|\Det\big(\alp(c_\eps)\big)\big| \ =\ e^{-\,\int_c\tet(h^{\Ea})+ \int_M\,\tet(h^{\Ea})\wedge X^*\Psi(g)}.
\end{equation}
Combining this equality with \refe{ration of metrics}
\eq{BZ}
	\frac{\|\cdot\|^\RS}{\|\cdot\|^{\M}_F}\ = \ 	e^{-\frac12\int_M\,\tet(h^{\Ea})\wedge X^*\Psi(g)},
\end{equation}
which is exactly the Bizmut-Zhang formula \cite[Theorem~0.2]{BisZh92}.

\def\cprime{$'$} \def\cprime{$'$} \newcommand{\noop}[1]{} \def\cprime{$'$}
\providecommand{\bysame}{\leavevmode\hbox to3em{\hrulefill}\thinspace}
\providecommand{\MR}{\relax\ifhmode\unskip\space\fi MR }
\providecommand{\MRhref}[2]{%
  \href{http://www.ams.org/mathscinet-getitem?mr=#1}{#2}
}
\providecommand{\href}[2]{#2}



\begin{thebibliography}{10}

\bibitem{BeGeVe}
N.~Berline, E.~Getzler, and M.~Vergne, \emph{Heat kernels and {Dirac}
  operators}, Springer-Verlag, 1992.

\bibitem{BisZh92}
J.-M. Bismut and W.~Zhang, \emph{An extension of a theorem by {Cheeger} and
  {M\"uller}}, Ast\'erisque \textbf{205} (1992).

\bibitem{BrKappelerRATshort}
M.~Braverman and T.~Kappeler, \emph{A refinement of the {Ray}-{Singer}
  torsion}, C.R. Acad. Sci. Paris \textbf{341} (2005), 497--502.

\bibitem{BrKappelerRATdetline_hol}
\bysame, \emph{Ray-{S}inger type theorem for the refined analytic torsion},
  \textbf{243} (2007), 232--256.

\bibitem{BrKappelerRATdetline}
\bysame, \emph{{Refined Analytic Torsion as an Element of the Determinant
  Line}}, Geometry \& Topology \textbf{11} (2007), 139–--213.

\bibitem{BrKappelerRAT}
\bysame, \emph{{Refined Analytic Torsion}}, J. Differential Geom. \textbf{78}
  (2008), no.~1, 193--267.

\bibitem{BurgheleaHaller_Euler}
D.~Burghelea and S.~Haller, \emph{{{Euler} Structures, the Variety of
  Representations and the {Milnor}-{Turaev} Torsion}}, Geom. Topol. \textbf{10}
  (2006), 1185–--1238.

\bibitem{Cheeger79}
J.~Cheeger, \emph{Analytic torsion and the heat equation}, Ann. of Math.
  \textbf{109} (1979), 259--300.

\bibitem{FarberTuraev00}
M.~Farber and V.~Turaev, \emph{Poincar\'e-{R}eidemeister metric, {E}uler
  structures, and torsion}, J. Reine Angew. Math. \textbf{520} (2000),
  195--225.

\bibitem{Franz35}
W.~Franz, \emph{{\"U}ber die torsion einer \"uberdeckung.}, Journal f\"ur die
  reine und angewandte Mathematik \textbf{173} (1935), 245--254.

\bibitem{Gilkey84}
P.~B. Gilkey, \emph{The eta invariant and secondary characteristic classes of
  locally flat bundles}, Algebraic and differential topology---global
  differential geometry, Teubner-Texte Math., vol.~70, Teubner, Leipzig, 1984,
  pp.~49--87.

\bibitem{GoldmanMillson88}
W.~Goldman and J.~Millson, \emph{The deformation theory of representations of
  fundamental groups of compact {K}\"ahler manifolds}, Inst. Hautes \'Etudes
  Sci. Publ. Math. (1988), no.~67, 43--96.

\bibitem{Ho04}
N.-K. Ho, \emph{The real locus of an involution map on the moduli space of flat
  connections on a {R}iemann surface}, Int. Math. Res. Not. (2004), no.~61,
  3263--3285.


\bibitem{KamberTondeur67}
F.~Kamber and Ph. Tondeur, \emph{Flat bundles and characteristic classes of
  group-representations}, Amer. J. Math. \textbf{89} (1967), 857--886.

\bibitem{MaZhang06eta}
X.~Ma and W.~Zhang, \emph{{$\eta$}-invariant and flat vector bundles}, Chinese
  Ann. Math. Ser. B \textbf{27} (2006), 67--72.

\bibitem{Muller78}
W.~M\"uller, \emph{Analytic torsion and {R}-torsion of {Riemannian} manifolds},
  Adv. in Math. \textbf{28} (1978), 233--305.

\bibitem{Muller93}
\bysame, \emph{Analytic torsion and {R}-torsion for unimodular representation},
  Jour. of AMS \textbf{6} (1993), 721--753.

\bibitem{Ponge-asymetry}
R.~Ponge, \emph{Spectral asymmetry, zeta functions, and the noncommutative
  residue}, Internat. J. Math. \textbf{17} (2006), 1065--1090.

\bibitem{Quillen85}
D.~Quillen, \emph{Determinants of {Cauchy}-{Riemann} operators over a {Riemann}
  surface}, Funct. Anal. Appl. \textbf{14} (1985), 31--34.

\bibitem{RaySinger71}
D.~B. Ray and I.~M. Singer, \emph{{R}-torsion and the {Laplacian} on
  {Riemannian} manifolds}, Adv. in Math. \textbf{7} (1971), 145--210.

\bibitem{Reidemeister35}
K.~Reidemeister, \emph{{H}omotopieringe und {L}insenr¬aume}, Hamburger Abhandl.
  \textbf{11} (1935), 102--109.

\bibitem{Seeley67}
R.~Seeley, \emph{Complex powers of an elliptic operator}, Proc. Symp. Pure and
  Appl. Math. AMS \textbf{10} (1967), 288--307.

\bibitem{ShabatSeveralVariables}
B.~V. Shabat, \emph{Introduction to complex analysis. {P}art {II}},
  Translations of Mathematical Monographs, vol. 110, American Mathematical
  Society, Providence, RI, 1992, Functions of several variables, Translated
  from the third (1985) Russian edition by J. S. Joel.

\bibitem{ShubinPDObook}
M.~A. Shubin, \emph{Pseudodifferential operators and spectral theory}, Springer
  Verlag, Berlin, New York, 1987.

\bibitem{Turaev86}
V.~G. Turaev, \emph{{R}eidemeister torsion in knot theory}, Russian Math.
  Survey \textbf{41} (1986), 119--182.

\bibitem{Turaev90}
\bysame, \emph{Euler structures, nonsingular vector fields, and
  {R}eidemeister-type torsions}, Math. USSR Izvestia \textbf{34} (1990),
  627--662.

\bibitem{Turaev01}
\bysame, \emph{Introduction to combinatorial torsions}, Lectures in Mathematics
  ETH Z\"urich, Birkh\"auser Verlag, Basel, 2001, Notes taken by Felix Schlenk.

\end{thebibliography}
\end{document}